\documentclass[11pt,a4paper]{article}
\usepackage[latin1]{inputenc}
\usepackage{amsmath}
\usepackage{amsfonts}
\usepackage{amssymb}
\author{S.J. Patterson}
\title{On a question of H.A. Schwarz\footnote{
The mathematics described here was presented and
discussed at the conference ``Thermodynamics Formalism - Applications
to Geometry, Number Theory and Stochastics'' (8th-12th July, 2019) at 
the Mittag-Leffler Institute, Djursholm.   The author thanks the 
institute and its staff for the very pleasant environment they provided.\newline
The work described here is part of a project to document the life
and work of Kurt Heegner whose papers have been deposited in the
Handschriftenabteilung of the SUB Göttingen.  This was a collaboration 
with Hans Opolka and Norbert Schappacher; H. Opolka undertook the task of
producing transcriptions of Heegner's most important mathematical
manuscripts and gave the first analysis of them.  We shall 
not touch on the biography of Heegner here.\newline
The text includes two quotations from a manuscript in the 
Handschriftenabteilung of the SUB Göttingen.   I would like to thank
the staff there for their courteous help with this project. \newline
Two last words of thanks - first to Tim Browning (Bristol), who made 
some stimulating and helpful remarks. \newline  
Finally to the team behind PARI/gp which I have used extensively in the 
work described here.\newline
\textit{Added by the revision:} I would like to thank the referee for a very helpful report.}}

\date{}

\begin{document}

\maketitle
\begin{center}
\textit{To Manfred Denker on the occasion of his 75\textsuperscript{th}
birthday}
\end{center}

\section{Introduction}
The question to which the title of this paper refers is not usually 
associated with the name of Hermann Amandus Schwarz and this needs 
some explanation.   To explain the problem let us consider a 
tetrahedron with vertices $1,2,3,4$; let $d_{ij}$ (with $i<j$) be 
the length of the side $ij$.   Let $V$ denote the volume of the 
tetrahedron.  Let $CM_0(x_{12},x_{13},x_{14},x_{23},x_{24},x_{34})$ 
denote the determinant of the matrix
\[
\left[\begin{array}{ccccc}
0&1&1&1&1\\
1&0&x_{12}^2&x_{13}^2&x_{14}^2\\
1&x_{12}^2 &0&x_{23}^2&x_{24}^2\\
1&x_{13}^2 &x_{23}^2&0&x_{34}^2\\
1&x_{14}^2&x_{24}^2&x_{34}^2&0
\end{array}
\right].
\]
This turns out to be a sum of $22$ monomials of degree 6 the coefficients 
of which are $\pm2$; it is therefore convenient to write
\[
CM(x_{12},x_{13},x_{14},x_{23},x_{24},x_{34})=
CM_0(x_{12},x_{13},x_{14},x_{23},x_{24},x_{34})/2 .
\]
The Cayley-Menger formula, which is actually due to Lagrange (1773), 
states that
\[
(12V)^2 = CM(d_{12},d_{13},d_{14},d_{23},d_{24},d_{34}).
\]
The question referred to in the title is whether the variety given 
by 
\[
y^2=CM(x_{12},x_{13},x_{14},x_{23},x_{24},x_{34})
\]
is rational over $\mathbb Q$, that is, whether we can find 
six elements $\xi_1,\xi_2,\xi_3,\xi_4,\xi_5,\xi_6$ of 
\[ \mathbb Q(x_{12},x_{13},x_{14},x_{23},x_{24},x_{34})[y]/
(y^2-CM(x_{12},x_{13},x_{14},x_{23},x_{24},x_{34}))\]
so that this 
field is $\mathbb Q(\xi_1,\xi_2,\xi_3,\xi_4,\xi_5,\xi_6)$.
Note that \[CM(x_{12},x_{13},x_{14},x_{23},x_{24},x_{34})\] 
is irreducible, \cite{DAS}.   

The question as to the rationality of a variety is a subtle one.
Lüroth's theorem shows that in the case of one variable any subfield
of $\mathbb Q (x)$ of finite index is itself of the form $\mathbb Q(\xi)$.
In higher dimensions the analogue of Lüroth's theorem does not hold and 
the question, even over $\mathbb C$, is delicate and has received 
considerable attention in recent times; for a discussion of such problems
for classes of varieties not all that different from the one with which 
we shall be concerned see \cite{HT}.   As it happens our case is 
rather special and the central problem is the question as to the 
set of algebraic number-fields for which the variety is rational.

Before going further it is instructive to look at the case of triangles.
The analogue of the Cayley-Menger determinant is the Heron function,
\[
H(x_{12},x_{23},x_{13})=-
\left|\begin{array}{cccc}
0&1&1&1\\
1&0&x_{12}^2&x_{13}^2\\
1&x_{12}^2 &0&x_{23}^2\\
1&x_{13}^2 &x_{23}^2&0\\
\end{array}
\right|
\]
which is equal to $2 x_{12}^2 x_{23}^2+2 x_{23}^2 x_{13}^2+2 x_{13}^2
x_{12}^2-x_{12}^4-x_{23}^4-x_{13}^4$, or, if one prefers, to
$(x_{12}+x_{23}+x_{13})(x_{12}+x_{23}-x_{13})(x_{12}-x_{23}+x_{13})
(-x_{12}+x_{23}+x_{13})$.
If we have a triangle $1,2,3$ with side-lengths $d_{12},d_{23},d_{13}$
and area $A$ then the formula usually ascribed to Heron (but probably 
older) asserts $(4A)^2 =H(d_{12},d_{23},d_{13})$.

Note that $H$ is reducible although neither $CM$ nor its higher dimensional 
analogues are -- see \cite{DAS}.   In this case the corresponding variety 
is rational; the problem has a long history and attracted even the 
attention of Gauss who wrote about it in a letter to H. Schumacher dated
21st October, 1847 and gave a parametric solution.  This is essentially 
the same as one given by Schottky in \cite{FS2} who discusses it at some
length.   One should note it is far from true that any solution of either
diophantine problem corresponds to to a geometric solution.  In the
case of triangles with rational area (Heron triangles) we obtain
the geometric solutions by restricting the triplet $(d_{12},d_{23},d_{13})$
so that all components are positive and the three triangle inequalities are
satisfied, that is, six inequalities in all.

With these ideas we can now discuss Schwarz' role in the history
of this problem.   He moved to Berlin as Weierstrass' successor in
1892, just before he turned 50.   Two years earlier he published his 
Collected Works and published only one paper during his tenure as 
professor in Berlin.   He seems to have decided on playing the 
``Grand Old Man''.   He had married, in 1868, Kummer's daughter, 
Marie, and apparently felt himself obliged to uphold the Kummer
tradition.  He was very active in the Berlin Mathematische 
Gesellschaft and sought out promising young students and encouraged 
them.   Kurt Heegner was one of the last, probably the last, and wrote
warmly of Schwarz' efforts on his behalf.   Kummer had taken up the
problem of finding all quadrilaterals with rational sides and rational
diagonals.   This was a problem going back to the 7\textsuperscript{th}
century Indian mathematician Brahmagupta. Kummer \cite{EEK} introduced 
a novel technique, employing, as we now say, elliptic curves, to study 
it.   One can, apart from the problem of understanding which solutions 
are geometric, consider this as the question of finding all the rational
solutions of $CM(x_{12},x_{13},x_{14},x_{23},x_{24},x_{34})=0$.  Kummer 
was quite possibly aware of the ``Cayley-Menger'' formula in later years.
His work on Kummer surfaces has its origins in the theory of the Fresnel
wave surface and Hamilton's discovery of a singularity of this surface. 
From this he made in 1833 the prediction, spectacular at the time, of the 
phenomenon of conical refraction.   Much of Cayley's  work of the 1840s 
deals with tetrahedroids which embody the algebraic approach to the
wave surface.  Kummer put forward, in the Berlin Mathematische Gesellschaft,
the problem finding all rational tetrahedra, that is, again leaving aside
the problem of determining which solutions are geometric, of finding the
rational (or integral) solutions of 
$y^2=CM(x_{12},x_{13},x_{14},x_{23},x_{24},x_{34})$.
This became a favourite problem of Schwarz and he encouraged several young
mathematicians to work on it, including Kurt Heegner.   One of the people
who took up the gauntlet was Otto Schulz whose extremely informative
doctoral thesis \cite{OS} contains, in particular, a historical summary 
which is the source for the information given above.   The problem 
of the the rationality of the variety
$y^2=CM(x_{12},x_{13},x_{14},x_{23},x_{24},x_{34})$
is attributed by Schulz to Schwarz and one can presume that 
Schwarz was motivated by the Heron example.  Schottky, who was later 
also in Berlin, discussed, in 1916, the Heron example in connection 
with the approach we shall describe here \cite[p.150ff.]{FS2}.  Another
Berlin student, a contemporary of Schulz, was Fritz Neiß (Neiss), whose
much shorter thesis \cite{FN}, deals with the same problem and offers some 
interesting insights.  Schulz give examples of rational subvarieties of 
the variety above.  Indeed one class of examples goes back to an
extension, by Kummer, of a theorem of Brahmagupta which can be interpreted 
as asserting that the subvariety of 
$\{ CM(x_{12},x_{13},x_{14},x_{23},x_{24},x_{34})=0 \}$
defined by the quadrilateral being inscribed in a given circle is rational.
This can be extended to tetrahedra by an observation of Friedrich Ankum
which we shall describe presently.

Schulz actually proves \newline\noindent
\textbf{Theorem}(O. Schulz) \textit{The variety}
\[ y^2=CM(x_{12},x_{13},x_{14},x_{23},x_{24},x_{34})\] 
\textit{is rational over} $\mathbb Q(\sqrt{-1})$.  
\textit{Moreover} 
\[ y^2=-CM(x_{12},x_{13},x_{14},x_{23},x_{24},x_{34})\]
\textit{is rational over} $\mathbb Q$.

He was clearly dissatisfied with this result (as it involved
$\mathbb Q(\sqrt{-1})$ rather than $\mathbb Q$) and relegated it to 
an appendix, \cite[p.268ff]{OS}.  He did not include the last statement 
but it follows from his argument and clearly implies the first one.  
His argument is fairly simple and we shall come to the essential point 
later.  The purpose of this paper is to give a more sophisticated proof 
based on ideas of Heegner which are preserved in the Handschriftenabteilung
of the SUB Göttingen, Codex Ms. K. Heegner 1:32, 1:33 
``Rationale Vierecke und Tetraeder und ihre Beziehung zu den 
elliptischen und hyperelliptischen Funktionen'').  The first 
of these is a handwritten manuscript, the second  a typescript.
Heegner's paper illuminates the problem further by pointing out
the relevance of a theorem of F. Schottky on Weddle surfaces,
\cite{FS1,FS2} which he develops further.  Heegner's work
is not dated but at the end of Ms. K. Heegner 1:33 we find the 
following remarks: \newline
\textit{Das vorliegende Manuskript habe ich im Jahre 
1940 abgeschlossen.   Die Fortführung der Arbeit hat mich zu 
algebraischen Untersuchungen über Modulfunktionen mehrerer 
Veränderlichen angeregt.  Soweit diese Untersuchungen mit den 
rationalen Vierecken und Tetraedern in Beziehung stehen, soll 
sich auch der zweite Teil der Arbeit damit befassen.}\footnote{
I completed the the present manuscript in the year 1940.   The
continuation of the work animated me to algebraic investigations
concerning modular functions in several variables.   That part of
the investigations connected with rational tetrahedra and 
quadrilaterals is to be taken up in the second part of this work.}   
\newline 
From the text we surmise that some of the ideas go back to his student
days (1915).  He writes in the introduction: \newline 
\textit{Ich selbst fiel H.A. Schwarz in den Kriegsjahren 1915/16 dadurch 
auf, daß ich bei der Behandlung dieser Aufgaben geeignete rechtwinklige 
Koordinaten einführte. Diese Untersuchungen haben mir mannigfache Anregung 
zu anderen Arbeiten gegeben. In größeren Zeitabständen bin ich jedoch 
zu diesen Aufgaben immer wieder zurückgekehrt und habe sie mit den
verschiedensten geometrischen Tatsachen in Beziehung gebracht.}
\footnote{I myself attracted the attention of H.A. Schwarz in the
war years 1915/1916 through treating this problem with appropriate
rectangular coordinates.   These investigations stimulated me 
to further works.  I have returned to this problem over longer 
periods of time and have found connections with several geometric 
facts.} \newline
The second part of ``Rationale Vierecke...'' is not in the 
Nachlass.   The remarks at the end of \cite{KH} are a brief and rather 
cryptic indication of what Heegner had in mind.

The results of the Berlin group were not published in any mathematical
journals; perhaps had the First World War not broken out at this 
time this would have changed.   The theses of Schulz and Neiss were 
circulated in printed versions, as was required at the time, but
they escaped the notice of L.E. Dickson and his team preparing
the ``History of the Theory of Numbers'', \cite{DTN} and 
remained unknown outside a very small circle of readers.  Heegner
never completed his investigations to his own satisfaction and poverty
and failing health in his later years meant that what he did achieve 
remained unpublished -- although he did drop a hint at the end of his 
final paper \cite{KH}\footnote{Interestingly this paper is included in 
the bibliography of \cite{CF} but is not cited in the body of the text 
as Prof. Flynn has kindly checked. Presumably the relevant passage 
was not included in the final text of the book.}.   His ideas remain
extremely interesting but now have to be embedded in the developments in
diophantine analysis of recent decades. One could see the present paper 
as ``Tales of forgotten genius'', to borrow a phrase from Ian Stewart,
\cite{B}.  There have been many papers over the years dealing with 
problems of rational quadrilaterals and tetrahedra but few have found 
the direction that was explored in Berlin.   Comparatively recently, in 
connection with Heron triangles, a related method was found 
by Robin Hartshorne and Ronald van Luijk, at that time a student of 
Hendrik Lenstra,  which brought the ideas from algebraic geometry into the 
study of this apparently fairly elementary problem - see \cite{HvL, vL}.

\section{Ankum's observation}

One important discovery of the group around Schwarz is due to
Friedrich Ankum; he noted that 
$CM(d_{12},d_{13},d_{14}+t,d_{23},d_{24}+t,d_{34}+t)$
is a quadratic function of $t$. Ankum never completed a thesis 
and his observation has been reported by Schulz, Neiss and Heegner.
Heegner recognised that the equation $CM(d_{12},d_{13},x,d_{23},y,z)=0$ 
in $(x,y,z)$ is a special type of Kummer surface, a tetrahedroid, and 
that Ankum's observation is an expression of the fact that the point 
at infinity  $[0:1:1:1]$ (with the usual convention of introducing 
the additional projective coordinate as the zeroth) is a node of the
tetrahedroid.   In fact tetrahedroids were the first class of Kummer
surfaces to be investigated, by Cayley, from 1846 onwards.   They are 
the algebraic version of the Fresnel wave surface when questions of 
the real locus are left out of the discussion.  For all of this see
\cite[Chapters IX,X]{H}.

The theory of tetrahedroids is not well covered in the contemporary 
literature.   There is a brief but useful account by A.-S. Elsenhans 
and J. Jahnel in \cite{EJ}.     There is a more extensive account 
by J.W.S. Cassels and E.V. Flynn in \cite{CF}, from a rather different 
point of view.  Elsenhans and Jahnel are interested in the arithmetic
of Kummer surfaces and tetrahedroids provide interesting examples. 
Cassels and Flynn's goal is to study the arithmetic of curves of
genus $2$.  Kummer surfaces are closely associated with them and 
in particular with their Jacobians but are more amenable to study.
Tetrahedroids are associated with curves with reducible Jacobians 
which happens quite often\footnote{There are two useful databases of 
curves of genus 2, one due to M. Stoll at \newline
\texttt{http://www.mathe2.uni-bayreuth.de/stoll}   \newline
and one due to A.R. Booker, J. Sijsling, A.V. Sutherland, J. Voight,
R. v. Bommel and D. Yasaki at the data base of L-Functions and modular 
forms, (see \cite{B}), at \newline
\texttt{http://www.lmfdb.org}. 
\newline   A perusal of these tables is most instructive.}.
Both of these accounts stress the analysis of a selected node and this 
proves not to be so convenient for our purposes.  Rather it turns 
out that the, more global, combinatorial approach taken by Hudson, 
\cite{H}, is appropriate.   For a more recent account see \cite{MGD}. 
We shall also need Cayley's theory of the symmetroid and the associated Weddle 
surface, briefly described in \cite{CF} but we need the account given in \cite{H}
or in the final chapter of \cite{J}.  Finally the two papers of W. 
Edge, \cite{WE1,WE2} provide an interesting perspective. The geometrical
literature does not pay much attention to fields of definition and
our account is a matter of making all of the details
concrete with the help of computer algebra, in the spirit
of \cite{CF}.

Let  
\[
T_{a,b,c}(X_0,X_1,X_2,X_3)=\frac 12 
\left|\begin{array}{ccccc}
0&1&1&1&X_0^2\\
1&0&a^2&b^2&X_1^2\\
1&a^2 &0&c^2&X_2^2\\
1&b^2 &c^2&0&X_3^2\\
X_0^2&X_1^2&X_2^2&X_3^2&0
\end{array}
\right|,
\]
i.e. $X_0^4 CM(a,b,X_1/X_0,c,X_2/X_0,X_3/X_0)$.  
Then $T_{a,b,c}(X_0,X_1,X_2,X_3)=0$ is the projective equation of 
a tetrahedroid and which we denote by $\mathcal T_{a,b,c}$; 
this is the family considered by Heegner.  The strategy we shall
now follow is the following. We consider our variety as fibered by 
the  $\mathcal T_{a,b,c}$ over the set of $(a,b,c)$.   We shall show 
that each $\mathcal T_{a,b,c}$ is rational in the appropriate sense and
that the rational functions realizing this are rational in $(a,b,c)$. 
The polynomial $T_{a,b,c}$ is a sum of $10$ monomials in $X_0,X_1,X_2,X_3$.

This yields somewhat more.  It is not difficult to see that the variety 
in $(a,b,c)$ determined by requiring that $H(a,b,c)$ lie in a 
fixed square class is rational (the argument is recalled at the end of 
Section 4) and it follows that the corresponding subvariety of 
$y^2=CM(x_{12},x_{13},x_{14},x_{23},x_{24},x_{34})$ is also rational.   
This will play no role here but it is relevant to the the idea Heegner 
refers to in the second quotation from his Nachlass (the introduction
of rectangular coordinates) given in the Introduction.

We shall work over a field $k$ of characteristic $0$ (usually 
$\mathbb Q$); we shall associate with the ``triangle'' $1 2 3$ the 
element $\delta_{a,b,c}=H(a,b,c) k^{\times 2}$  of 
$k^{\times}/k^{\times 2}$.  Neiss observed the following identity 
which he deduced from Sylvester's identity for determinants
\[
H(d_{12},d_{13},d_{23})H(d_{12},d_{24},d_{14})=D_{12}^2+
(2 d_{12})^2 CM(d_{12},d_{13},d_{14},d_{23},d_{24},d_{34}) 
\]
where $D_{12}=D_{12}(d_{12},d_{13},d_{14},d_{23},d_{24},d_{34})$ is 
\[
\left|\begin{array}{cccc}
0&1&1&1\\
1&0&d_{12}^2&d_{14}^2\\
1&d_{12}^2&0&d_{24}^2\\
1&d_{13}^2 &d_{23}^2&d_{34}^2
\end{array}
\right|.
\]
From this one can make two deductions.   
If $CM(d_{12},d_{13},d_{14},d_{23},d_{24},d_{34})=0$, which would 
be the case if we were dealing with a plane quadrilateral,
then the square-classes $\delta_{a,b,c}$ for all three element
subsets of $\{1,2,3,4\}$ are equal.   If
$CM(d_{12},d_{13},d_{14},d_{23},d_{24},d_{34})\neq 0$
then square classes $\delta_{a,b,c}$ of two triangles are equal 
if $-1$ is a square in $k$ and differ multiplicatively by a norm 
from $k(\sqrt{-1})^{\times}$ otherwise. Here we  assume that none of 
the corresponding $H$ are zero. 
 
Finally Schulz noted that 
$D_{12}(d_{12},d_{13},d_{14}+s,d_{23},d_{24}+s,d_{34}+s)$
is of the form $As+B$ where 
\[
A=-2((d_{12}^2-d_{13}^2+d_{23}^2)d_{14}^2+(d_{12}^2+d_{13}^2-d_{23}^2)d_{24}^2
-2d_{12}^2d_{34}^2),
\]
a refined version of Ankum's observation.   This shows that 
$CM(d_{12},d_{13},d_{14}+s,d_{23},d_{24}+s,d_{34}+s)$ is 
a polynomial of degree $2$ in $s$ and the leading coefficient is 
positive.   From this and the theory of binary quadratic forms one 
can derive, from any given rational solution an infinite family of 
rational solutions of 
$y^2=CM(x_{12},x_{13},x_{14},x_{23},x_{24},x_{34})$.

\section{Some classical mathematics}
In this section we shall recall some classical ideas, mainly due to 
Cayley, which we shall use.   Let $\mathcal P$ be a set of 
six points in $\mathbb P^3$, no four of which lie in a plane.   Then 
the space of quadratic forms which vanish at these points is 
four dimensional; let $S_1,S_2,S_3,S_4$ be a basis.   If we take a 
point $\xi \in \mathbb P^3$ outside $\mathcal P$ then the subspace 
of quadratic forms vanishing at $\xi$ is one-dimensional.   However, by 
Bézout's theorem, \cite[\S I.7]{HAG}, two quadrics in $\mathbb P^3$ meet
in eight points.   Thus there is a further point $\xi'$ associated with $\xi$.   
The map 
$\mathbf S:\mathbb P^3 \setminus \mathcal P\rightarrow \mathbb P^3;
\mathbf x \mapsto  [S_1(\mathbf x):S_2(\mathbf x):
S_3(\mathbf x):S_4(\mathbf x)] $ is therefore of degree two.   It 
corresponds to a field extension which will then be quadratic. As 
such it is is Galois and associated with a (non-linear) involution of 
$\mathbb P^3$. The details can be found in \cite{Sn}.   The
explicit formul\ae\ for the involution are rather complicated and
seem only to be accessible through the use of computer algebra.   It
is also questionable whether they are of much use.

We recall that a \textit{trope} of a Kummer surface is the 
intersection of the Kummer surface with a plane tangential to it 
along a conic.   There are sixteen tropes and they correspond to 
the nodes of the dual Kummer surface -  \cite[Ch.4]{CF}.

To this framework one can associate three varieties, $K_\mathcal P, 
K_\mathcal P^*$ and $W_\mathcal P$.   They are defined as follows.  
$K_\mathcal P$ is the vanishing set of the degree four polynomial
$\det(z_1 S_1 + z_2 S_2+ z_3 S_3 +z_4 S_4)$ where the $S_i$ are
considered as matrices.   This is called the \textit{symmetroid}.
It is a Kummer surface and fairly easy to determine in terms of 
$\mathcal P$.   Let $JS(x)=\det((\partial S_i/\partial z_j)_{1 \le i,j \le 4} )$. 
Then, by definition, the Weddle surface $W_\mathcal P$ is the vanishing 
locus of $JS(x)$.  There is a 1--1 correspondence between $K_\mathcal P$
and $W_\mathcal P$, see \cite[(5.1.5)]{CF}; it is cubic and arises from 
the minors of the determinant.   Finally let $K_\mathcal P^*$ be the 
image of $W_\mathcal P$ under $\mathbf S$.   It is 
again a Kummer surface, the dual of $K_\mathcal P$, and the map 
$\mathbf S$ is a partial desingularization of $K_\mathcal P^*$; fifteen 
of the nodes are blown up to lines in $W_\mathcal P$.   The map can be described as 
follows:  for each pair of points $\{P,P'\}$ in $\mathcal P$ the image
under $\mathbf S$ of the line joining $P$ and $P'$ is a single
point, and is then a node of  $K_\mathcal P^*$.  This yields the fifteen
nodes; the final one is not in the image of $\mathbf S$. For any partition 
of $\mathcal P$ into two disjoint subsets $\Pi,\Pi'$, each of three elements, 
the planes defined by $\Pi$ and $\Pi'$ meet in a line and the image of 
this line is a trope  in $K_\mathcal P^*$.  We obtain ten of the sixteen
tropes in this way.  The complete picture, which gives a beautiful geometric
interpretation of the configuration is given, rather tersely in
\cite[pp. 166, 167]{H}, especially the table at the top of p.167.  The memoir 
\cite{MGD} deals rather with complete desingularizations of Kummer surfaces and 
does not cover the theory of Weddle surfaces. 

The combinatorics of these nodes and tropes describe the relations between 
them in a very convenient fashion.   Tetrahedroids are characterised amongst 
Kummer surfaces by the fact that the nodes split into four sets of four, 
each of which is contained in a plane.   The four planes define the 
tetrahedron giving the tetrahedroid its name.  

In our case the combinatorics can be described explicitly by the following 
two tables.   They give respectively the nodes and the tropes of the 
tetrahedroid. 
\[
\tiny
\begin{array}{||ccc||ccc||}
\hline \hline
&&&&&\\
\emptyset \cup \{1,2,3,4,5,6\} &\longleftrightarrow& [0:1:1:1]&
\{4,6\}\cup \{1,2,3,5\} & \longleftrightarrow &[1:0:a:b]\\
\{3,4\}\cup \{1,2,5,6\} &\longleftrightarrow &[0:1:-1:1]&
\{4,5\}\cup \{1,2,3,6\} &\longleftrightarrow &[1:0:-a:b]\\
\{5,6\}\cup \{1,2,3,4\} &\longleftrightarrow &[0:1:1:-1]&
\{3,6\}\cup \{1,2,4,5\} &\longleftrightarrow &[1:0:a:-b]\\
\{1,2\}\cup \{3,4,5,6\} &\longleftrightarrow &[0:1:-1:-1]&
\{3,5\}\cup \{1,2,4,6\} &\longleftrightarrow &[1:0:-a:-b]\\
&&&&&\\
\hline
&&&&&\\
\{2,5\}\cup \{1,2,3,5\} &\longleftrightarrow& [1:a:0:c]&
\{1,3\}\cup \{1,2,3,5\} &\longleftrightarrow &[1:b:c:0]\\
\{2,6\}\cup \{1,2,5,6\} &\longleftrightarrow &[1:-a:0:c]&
\{1,4\}\cup \{1,2,3,6\} &\longleftrightarrow &[1:-b:c:0]\\
\{1,5\}\cup \{1,2,3,4\} &\longleftrightarrow &[1:a:0:-c]&
\{2,3\}\cup \{1,2,4,5\} &\longleftrightarrow &[1:b:0:-c]\\
\{1,6\}\cup \{3,4,5,6\} &\longleftrightarrow &[1:-a:0:-c]&
\{2,4\}\cup \{1,2,4,6\} &\longleftrightarrow &[1:-b:0:-c]\\ 
&&&&&\\
\hline \hline
\end{array}\ ,
\]
\[
\tiny
\begin{array}{||ccc||ccc||}
\hline \hline
&&&&&\\
\{1,3,5\}\cup \{2,4,6\} &\longleftrightarrow &[0:c:b:a]&
\{1,4,6\}\cup \{2,3,5\} &\longleftrightarrow &[c:0:1:1]\\
\{2,3,5\}\cup \{1,4,6\} &\longleftrightarrow &[0:c:-b:a]&
\{1\}\cup \{2,3,4,5,6\} &\longleftrightarrow &[c:0:-1:1]\\
\{1,3,6\}\cup \{2,4,5\} &\longleftrightarrow &[0:c:b:-a]&
\{2\}\cup \{1,3,4,5,6\}&\longleftrightarrow &[c:0:1:-1]\\
\{1,4,5\}\cup \{1,3,6\} & \longleftrightarrow &[0:c:-b:-a]&
\{1,3,4\}\cup \{2,5,6\} &\longleftrightarrow &[c:0:-1:-1]\\
&&&&&\\
\hline
&&&&&\\
\{1,2,4\}\cup \{3,5,6\} &\longleftrightarrow& [b:1:0:1]&
\{1,2,6\}\cup \{3,4,5\} &\longleftrightarrow &[a:1:1:0]\\
\{3\}\cup \{1,2,4,5,6\} &\longleftrightarrow &[b:-1:0:1]&
\{5\}\cup \{1,2,3,4,6\} &\longleftrightarrow &[a:-1:1:0]\\
\{4\}\cup \{1,2,3,5,6\} &\longleftrightarrow &[b:1:0:-1]&
\{6\}\cup \{1,2,3,4,5\} &\longleftrightarrow &[a:1:-1:0]\\
\{1,2,3\}\cup \{4,5,6\} &\longleftrightarrow &[b:-1:0:-1]&
\{1,2,5\}\cup \{3,4,6\} &\longleftrightarrow &[a:-1:-1:0]\\ 
&&&&&\\
\hline \hline
\end{array}\ .
\]
The four groups of four belong to the faces of the two tetrahedra.
It is straightforward from this to write down the incidence matrix
explicitly.  It is valid for all Kummer surfaces, not just 
tetrahedroids.   

The locus of the fixed points of the involution described above 
are contained in $W_\mathcal P$ and, as it is irreducible, is
identical with it.   An argument involving Galois theory and 
Bézout's theorem shows that $JS(x)^2$ is in the ring generated
by $S_1, S_2, S_3, S_4$.   That is, there is a polynomial $F$, which is
of degree four and with coefficients in $\mathbb Q(a,b,c)$, so that 
$JS(x)^2=F(S(x))$  -  see \cite[p.170]{H}, \cite{FS1,FS2}.   This result 
seems to be due to Schottky and we shall refer to it as ``Schottky's theorem''.
Whereas the nodes of $K_\mathcal P^*$ are determined by the constructions 
above the polynomial $F$ is not so easily determined in general.  However, in 
our case it is quite easy as the dual $\mathcal T_{a,b,c}^*$  of 
$\mathcal T_{a,b,c}$ turns out to be the zero set of 
\[
T_{a,b,c}^*(X_0,X_1,X_2,X_3)=\frac 12 
\begin{array}{|ccccc|}
0&c^2&b^2&a^2&X_0^2\\
c^2&0&1&1&X_1^2\\
b^2&1&0&1&X_2^2\\
a^2&1 &1&0&X_3^2\\
X_0^2&X_1^2&X_2^2&X_3^2&0
\end{array}\ .
\]
We note that 
\[T_{a,b,c}^*(abcX_0,cX_1,bX_2,aX_3)=(abc)^2T_{a,b,c}(X_0,X_1,X_2,X_3),\]
that 
\[T_{\lambda a,\lambda b,\lambda c}(X_0,X_1,X_2,X_3)=
\lambda^2 T_{a,b,c}(\lambda X_0,X_1,X_2,X_3)
\]
and that
\[T_{\lambda a,\lambda b,\lambda c}^*(X_0,X_1,X_2,X_3)=
\lambda^2 T_{a,b,c}^*(\lambda^{-1} X_0,X_1,X_2,X_3).
\]
The classical theory (\cite[p.40 and passim]{CF}, \cite[p.170]{H}) shows 
that $K_\mathcal P^*$ is the dual of $K_\mathcal P$  They are 
isomorphic over the complex numbers but not necessarily otherwise,
\cite[Ch. 4, esp. Theorem 4.5.1]{CF}.   In the case of 
tetrahedroids they actually are isomorphic over the field of definition.  
The formul\ae\ above demonstrate this.

We asserted above that $K_\mathcal P$ can be determined from 
$\mathcal P$.   To do this we can arrange the coordinate system
so that four of the points of $\mathcal P$ are $[1:0:0:0]$,
$[0:1:0:0]$, $[0:0:1:0]$ and $[0:0:0:1]$.   There remain two 
additional points which we write as $[p_1:p_2:p_3:p_4]$ and 
$[q_1:q_2:q_3:q_4]$.   In fact, by assumption, none of the 
$p_j$ are zero and so we could also arrange that
$[p_1:p_2:p_3:p_4]=[1:1:1:1]$; this makes computations easier but
sometimes obscures the structure of the formul\ae .  In the general 
case we can (assuming $p_3q_2 \neq p_2q_3$) take, as a basis the set 
of quadratic forms
\[
\small
\begin{array}{l}
S_1=(p_4q_3-p_3q_4)x_1x_2 +(p_2q_4-p_4q_2)x_1x_3+(p_3q_2-p_2q_3)x_1x_4,\\
S_2=p_3q_3(p_2q_1-p_1q_2)x_1x_2
+p_2q_2(p_1q_3-p_3q_1)x_1x_3+p_1q_1(p_3q_2-p_2q_3)x_2x_3,\\
S_3=(p_2q_4q_1q_3-p_1p_3q_2q_4)x_1x_2
+p_2q_2(p_1q_4-p_4q_1)x_1x_3+p_1q_1(p_3q_2-p_3q_3)x_2x_4,\\
S_4=p_3q_3(p_4q_1-p_1q_4)x_1x_2
+(p_1p_2q_3q_4-p_3p_4q_1q_2)x_1x_3+q_1(q_2-q_3)x_3x_4.
\end{array}
\]
One can evaluate $JS(\mathbf x)$ from this but one knows by  \cite[p.170]{H} that it is, 
up to a multiple independent of $\mathbf x$,
\[
W_{\mathbf p,\mathbf q}(\mathbf x)=
\begin{array}{|cccc|}
x_1^2&p_1x_1&q_1x_1&p_1q_1\\
x_2^2&p_2x_2&q_2x_2&p_2q_2\\
x_3^2&p_3x_3&q_3x_3&p_3q_3\\
x_4^2&p_4x_4&q_4x_4&p_4q_4
\end{array}\ .
\]

The determination of $K_\mathcal P$ is as follows.   We 
let $V=(v_{i,j})_{1 \le i,j \le 4}$ be a symmetric matrix with 
$v_{i,i}=0$ $1\le i \le 4$. Then $_{\mathcal P}$ can be identified
with the intersection of the hypersurface in $\mathbb P^5$ given 
by $\det(V)=0$ with the two hyperplanes $\sum_{i<j}p_ip_j v_{i,j} =0$ 
and $\sum_{i<j}q_iq_j v_{i,j} =0$.   This is associated with 
the so-called ``irrational form'' (\cite[\S 19]{H}) of the equation 
of a Kummer surface.   Also from this one can give the map from 
$K_\mathcal P \rightarrow W_\mathcal P$ explicitly; for all
this see \cite[pp.171,172]{H}.   The maps are cubic, as noted above.

One can derive the ``irrational form'' relatively easily 
from information about the nodes and tropes of the Kummer 
surface - \cite[pp.34--36]{H}.   The ``irrational form''
involves three products of pairs of linear forms defining 
tropes.  Write these as $L_1(\mathbf x) L_1'(\mathbf x)$,
$L_2(\mathbf x) L_2'(\mathbf x)$ and
$L_3(\mathbf x) L_3'(\mathbf x)$; the equation takes the form
\[
\sqrt{L_1(\mathbf x) L_1'(\mathbf x)}+
\sqrt{L_2(\mathbf x) L_2'(\mathbf x)}+
\sqrt{L_3(\mathbf x) L_3'(\mathbf x)}=0.
\]
In a form not involving questionable irrationalities this becomes 
\[
(L_1(\mathbf x) L_1'(\mathbf x)+
L_2(\mathbf x) L_2'(\mathbf x)-
L_3(\mathbf x) L_3'(\mathbf x))^2-
4L_1(\mathbf x) L_1'(\mathbf x)L_2(\mathbf x) L_2'(\mathbf x)=0;
\]
compare this with the formula for $T_{a,b,c}$ (Schulz' identity) given 
below.

We shall take
\[
\begin{array}{cc}
L_1=-(a-b+c)(c x_0+x_2+x_3), & L_1'=(-a+b+c)(-c x_0+x_2+x_3)\\
L_2=(a+b+c)(c x_0+x_2-x_3), & L_2'=(a+b-c)(-c x_0+x_2-x_3)\\
L_3=2(-c x_1-b x_2+a x_3), & L_3'=2(c x_1-b x_2 +a x_3)
\end{array}
\]
One verifies that $L_1+L_1'+L_2+L_2'+L_3+L_3'=0$ and with \newline 
$q_1=(a+b-c)(a-b+c)$, $q_2=-(a+b+c)(-a+b+c)$,\newline $q_3=-(-a+b+c)(a+b-c)$
and $q_4=(a+b+c)(a-b+c)$ one has
\[
L_1q_2q_3+L_1'q_1q_4+L_2q_1q_3+L_2'q_2q_4+L_3q_1q_2+L_3'q_3q_4=0.
\]

One finds now
\[
16c^2T_{a,b,c}=-((L_1L_1'+L_2L_2'-L_3L_3')^2-4 L_1L_1'L_2L_2).
\]
The identity Schulz used is a variant of this. His version of Ankum's observation
follows from it.

It follows from \cite[p.171]{H} that, by a suitable choice of basis  
$\mathcal P =\{P_1,P_2,P_3,P_4,P_5,P_6\} $ corresponding
to $\mathcal T_{a,b,c}$ can be taken to be the four basis points, which
we write as above with $[p_1:p_2:p_3:p_4]=[1:1:1:1]$. One can now calculate 
the nodes of the Kummer surface which is the image of the Weddle surface
under $\mathbf S$ .   The four planes of the tetrahedron are simple to 
find and with respect to the natural set of coordinates associated with 
them we can identify the image as $\mathcal T_{2a,2b,2c}^*$.
All of the coordinates are rational in $a,b,c$.   The condition that 
the Kummer surface be a tetrahedroid is that $q_1q_2=q_3q_4$, see
\cite[p.366]{Sn}. Note that $q_1q_2=-H(a,b,c)$.  The calculations 
at this point are routine and the results are not so elegant that they 
need be exhibited here.  It may be the case that these are contained 
in the formul\ae\ in \cite[\S 11]{HFB} -- see \cite[p.40]{CF}
and \cite[p.953]{WE1}.  I have not attempted to verify this.

\section{Completion of the proof of Schulz' theorem}
The results of the previous section shows that if $\mathcal P$
and $\mathbf S= (S_1,S_2,S_3,S_4)$ are as above then there 
is an element $\gamma \in \mathrm{GL}(4,\mathbb Q(a,b,c))$ so 
that for some element $k \in \mathbb Q(a,b,c)^\times$ one 
has 
\[
W_{\mathbf 1,\mathbf q}(\mathbf x)^2 = 
k T_{a,b,c}(\gamma(\mathbf S(\mathbf x))),
\]
where $\mathbf 1=(1,1,1,1)$ and $\mathbf q =(q_1,q_2,q_3,q_4).$
We can make all this explicit. We introduce the system $\mathbf S^*$:
\[
\begin{array}{rcl}
S_1^*(\mathbf X) &=& 2(X_1X_2-X_3X_4),\\
S_2^*(\mathbf X) &=& (a+b+c)X_1X_3-(-a+b+c)X_1X_4-(a-b+c)X_2X_3\\
&&\hspace{7cm}-(a+b-c)X_2X_4,\\
S_3^*(\mathbf X) &=&-2aX_1X_2+(a+b+c)X_1X_3-(-a+b+c)X_1X_4\\
&&\hspace{2cm}+(a-b+c)X_2X_3+(a+b-c)X_2X_4-2aX_3X_4,\\
S_4^*(\mathbf X) &=&-2bX_1X_2+(a+b+c)X_1X_3+(-a+b+c)X_1X_4\\
&&\hspace{2cm}-(a-b+c)X_2X_3+(a+b-c)X_2X_4-2bX_3X_4. 
\end{array}
\]
Then by a judicious comparison of coefficients we find 
\[
c^2 T_{a,b,c}(S^*_1(\mathbf X),S^*_2(\mathbf X),
S^*_3(\mathbf X),S^*_4(\mathbf X))=
-W_{\mathbf 1,\mathbf q}(\mathbf X)^2.
\]
This formula is so simple that it is easy to verify directly.
As the quadratic function $\mathbf S^*$, depends rationally on $a,b,c$, 
as does $W_{\mathbf 1,\mathbf q}$. Schulz' theorem follows directly.

It may also be noted that this theory also allows one to 
make the inverse map explicit.   We shall return to this question 
in a future publication.  

The tetrahedroids have nodes and tropes which are rational (see 
\cite{EJ}) but there are other constructs which are not 
rational.   Heegner gave another representation of the 
varieties $CM(\mathbf d)=0$ and $V^2=CM(\mathbf d)$ involving
$\mathbb Q(\sqrt{-H(d_{12},d_{13},d_{23})})$ or one of the
analogues for the other faces.   In the case of rational 
quadrilaterals this is independent of the choice of ``face''; 
in the case of rational tetrahedra it is now the case that the factors 
corresponding to two different faces can differ multiplicatively 
by a norm from $\mathbb Q(\sqrt{-1})^\times$ (see \S 2).  Heegner stressed 
the significance of classifying the solutions by, in the first case, 
the quadratic field and, in the second case, by the associated 
quaternion algebra.  There are quite a number of still unanswered 
questions which arise in this context.   One is whether there is a 
combinatorial structure such as the theory of the Markoff tree in 
the case of the Markoff equation and its generalizations, \cite[Ch. 2]{C},
\cite[pp.106--110]{M}.  In the penultimate paragraph of \cite{KH} Heegner
indicates that, in the case of rational quadrilaterals, he has constructed
an algorithm, but exactly what this is remains unclear.   The  present author 
has not yet managed to identify it in his notes. 

Finally we note that the proof of Schulz' theorem given here is 
a geometrical version of the original proof.   It does not rule
out that either the variety considered in the theorem, or the 
more refined varieties considered by Heegner are rational over 
$\mathbb Q$.     In the case of Heron triangles both $Y^2=H(a,b,c)$ 
and $Y^2=-H(a,b,c)$ are rational.   To see this let 
$U=(-a+b+c)/(a+b+c)$, $V=(a-b+c)/(a+b+c)$ and $Z=Y/(a+b+c)^2$.   
The two equations become $Z^2=UV(2-U-V)$ and $Z^2=-UV(2-U-V)$.   
For fixed $V$ these are quadratic in $U$ and have rational points, 
$U=0$ and $U=2-V$.  They are then rational and with varying $V$ we 
get the required parametrization.   The question as to the rationality 
of $V^2=CM(\mathbf d)$ itself over $\mathbb Q$, or, for that matter,
over $\mathbb R$, is still very much open and the example of Heron 
triangles shows that it is delicate.   One notes that the real
locus of $T_{a,b,c}(\mathbf x)=0$ is connected; if one removes the 
nodes it becomes the union of eight four-fold punctured spheres.   
If it were not connected it would represent an obstruction to
$y^2=T_{a,b,c}(\mathbf x)$ being rational.  At the present the 
evidence, such as it is, is against a positive resolution of Schwarz'
problem over $\mathbb Q$, or even $\mathbb R.$

\medskip
\noindent
Mathematisches Institut\\
Bunsenstr. 3--5\\
37073 G\"ottingen\\
Germany

\noindent
e-mail:\texttt{spatter@gwdg.de}

\end{document}